\documentclass[12pt]{article}

\usepackage{a4wide}
\usepackage{amssymb}
\usepackage{amsfonts}
\usepackage{amsmath}
\input xy
\xyoption{arrow} \xyoption{matrix}

\date{}

\newtheorem{proposition}{Proposition}[section]
\newtheorem{theorem}[proposition]{Theorem}
\newtheorem{lemma}[proposition]{Lemma}

\newtheorem{corollary}[proposition]{Corollary}

\def\der{\partial }

\def\nFM0{{\nu }_{F,M_0}}
\def\nFN0{{\nu }_{F,N_0}}
\def\nGN0{{\nu }_{G,N_0}}

\def\N0{ {\bf N}_0 }
\def\R{\mathbb{R}}

\def\g{\gamma}

\def\ra{\rightarrow}

\def\Xpm{X^{\pm }}

\def\s{\sigma}
\def\Z{\mathbb{Z}}

\def\l1{{\lambda}_1}

\def\a{\alpha}
\def\a0{ {\alpha }_0}
\def\a1{ {\alpha }_1}

\def\l{\lambda}
\def\o{\omega}

\def\nFGM0{{\nu }_{F,G,M_0}}

\def\nFN0{{\nu}_{F,N_0}}

\def\sm{{\sigma}^m}

\def\sm1{{\sigma}^{-1}}

\def\smtp1{{\sigma}^{-t+1}}

\def\o{\omega }
\def\S1{S^{-1}}

\def\Xpm1{X^{\pm 1}_1}

\def\sPM1{{\sigma }^{\pm 1}}
\def\sMP1{{\sigma }^{\mp 1 }}

\def\d{\delta}

\def\di{{\rm d.ind}}

\def\L{\Lambda}

\def\G{\Gamma}

\def\Ytm1{Y^{t-1}}
\def\Yim1{Y^{i-1}}

\def\Aut{{\rm Aut}}

\def\D{ \Delta }

\def\SL2Z{ {\rm SL}_2({\bf Z}) }

\def\th{ \theta }

\def\Gp1{ G^{1 , 1 } }
\def\P11{ P^{-1 , 1 } }
\def\Pp1{ P^{1 , 1 } }

\def\th{\theta}

\def\nCLsr{{}^\nu\kern-2pt {\cal L}^{\sigma , \rho  }}
\def\nP{{}^\nu \kern-2pt P}
\def\nL{{}^\nu\kern-2pt L}
\def\nLL{{}^\nu\kern-2pt \Lambda}
\def\nPsr{{}^\nu\kern-2pt P^{\sigma , \rho  }}
\def\nLsr{{}^\nu\kern-2pt L^{\sigma , \rho  }}
\def\nuCL{{}^\nu\kern-2pt  {\cal L}}
\def\nCLsr{{}^\nu\kern-2pt {\cal L}^{\sigma , \rho  }}
\def\nCL1m{{}^\nu\kern-2pt {\cal L}^{-1 , 1  }}
\def\x1nu{x^\frac{1}{\nu}}
\def\xm1nu{x^{-\frac{1}{\nu}}}

\def\ra{\rightarrow }

\def\CP{ {\cal P}}

\def\nAM0{{\nu }_{{\cal A},M_0}}
\def\nAN0{{\nu }_{{\cal A},N_0}}

\def\CP{ {\cal P }}
\def\det{ {\rm det }}

\def\SL{{\rm SL}}

\def\di!{\frac{\der^i}{i!}}
\def\dik!{\frac{\der^k_i}{k!}}

\def\N{\mathbb{N}}

\def\0{\overline{0}}
\def\1{\overline{1}}

\def\Ln1{\L_{n,\overline{1}}}

\def\a1{a_{\overline{1}}}

\def\S{\Sigma}

\def\Aff{{\rm Aff}}

\def\Q{\mathbb{Q}}

\def\Mon{{\rm Mon}}
\def\Dec{{\rm Dec}}
\def\Sign{{\rm Sign}}
\def\sign{{\rm sign}}
\def\vn1{\overrightarrow{n-1}}
\def\CQ{{\cal Q}}

\begin{document}

\author{V. V. \  Bavula 
}

\title{Factorization of monomorphisms of a polynomial  algebra in one variable}

\maketitle
\begin{abstract}
Let $K[x]$ be a polynomial algebra in a variable $x$ over a
commutative $\Q$-algebra $K$, and $\G'$ be the monoid of
$K$-algebra monomorphisms of $K[x]$ of the type $\s : x\mapsto
x+\l_2x^2+\cdots +\l_nx^n$, $\l_i\in K$, $\l_n$ is a unit of $K$.
It is proved that for each $\s \in \G'$ there are only finitely
many distinct decompositions $\s = \s_1\cdots \s_s$ in $\G'$.
Moreover, each such a decomposition is uniquely determined by the
degrees of components: if $\s = \s_1\cdots \s_s= \tau _1\cdots
\tau_s$ then $\s_1=\tau_1, \ldots , \s_s=\tau_s$  iff $\deg (\s_1
)=\deg (\tau_1), \ldots , \deg (\s_s)=\deg (\tau_s)$. Explicit
formulae are given for the components $\s_i$ via the coefficients
$\l_j$ and the degrees $\deg (\s_k)$ (as an application of the
inversion formula for polynomial automorphisms in {\em several}
variables from \cite{Bav-inform}). In general, for a polynomial
there are no formulae (in radicals) for its divisors (elementary
Galois theory). Surprisingly, one can write such formulae where
instead of the product of polynomials one considers their
composition (as polynomial functions).



$${\bf Contents}$$
\begin{enumerate}
\item Introduction. \item Proof of Theorem \ref{27Dec06}.\item
Formulae for the components $\s$ and $\tau$ in $\d = \s\tau$.
\item Necessary conditions for irreducibility of a polynomial.
\end{enumerate}
\end{abstract}


\section{Introduction}
Throughout,  $K$ is a commutative $\Q$-algebra (if it is not
stated otherwise) with the group of units $K^*$, $K[x]$ is a
polynomial algebra over $K$ in a single variable $x$,
$\Aut_K(K[x])$ and $\Mon_K(K[x])$ are the group of automorphisms
and the monoid of monomorphisms of the polynomial algebra $K[x]$
respectively.
$$\G':=\G'(K):=\{\s \in  \Mon_K(K[x])\, | \,  \s :x\mapsto x+\l_2x^2+\cdots
+\l_nx^n, \l_i\in K, \l_n\in K^*\}$$ is the submonoid of
$\Mon_K(K[x])$, and $\deg (\s ) := \deg (\s (x))$ is called the
{\em degree} of $\s$. For $\s_1, \ldots ,\s_s\in \G'$,
\begin{equation}\label{degm}
\deg (\s_1\cdots \s_s) = \deg (\s_1)\cdots \deg (\s_s).
\end{equation}
The group $\Aut_K(K[x])$ contains the  {\em affine}  group $\Aff
:=\{ \s : x\mapsto ax+b\, | \, a\in K^*, b\in K\}$. If $K$ is a
field of characteristic zero then $\Aut_K(K[x])= \Aff$. If, in
addition, the field $K$ is algebraically closed then
$\Mon_K(K[x])= \Aff\times_{ex} \G'$ is the exact product of
monoids, i.e. each monomorphism $\s \in \Mon_K(K[x])$ is a unique
product $\s = \tau \g$ for some $\tau \in \Aff$ and $\g \in \G'$.

 The monoid $\G'$ is large, it is an infinite dimensional
 algebraic monoid. The submonoid  $M$ of $\G'$ generated by the
 monomorphisms $\{\s_{\l x^n} : x\mapsto x+\l x^n\, | \, \l\in K^*,
 n\geq 2\}$ is a {\em free} monoid (Theorem \ref{23Dec06}) with
 the free generators  $\{\s_{\l x^n}\}$.


For an element $\s \in \G'\backslash \{ e\}$ (where $e$ is the
identity of $\G'$), the set
$$ \Dec (\s ) := \{ (\s_1, \ldots , \s_s)\, | \, \s_1\cdots
\s_s=\s , s\geq 1 , \s_i\in \G'\backslash \{ e\}\}$$ is called the
{\em decomposition set} for $\s$,  the set $$\Sign (\s ) :=
\{(\deg (\s_1), \ldots , \deg(\s_s))\, | \, (\s_1, \ldots ,
\s_s)\in \Dec (\s)\}$$is called the {\em signature} of $\s$, and
the map
$$ \sign := \sign_\s : \Dec (\s ) \mapsto \Sign (\s ) , \;\; (\s_1, \ldots ,
\s_s)\mapsto (\deg (\s_1), \ldots , \deg(\s_s)),$$ is called the
{\em signature map}. It is obvious that the signature map is a
surjection and the signature  $\Sign (\s )$ is a {\em finite} set
since, for each $(\s_1, \ldots , \s_s)\in \Dec (\s )$,  $\deg (\s
) = \deg (\s_1)\cdots \deg(\s_s)$. The next theorem shows that the
signature map is a bijection, i.e. each decomposition $\s =
\s_1\cdots \s_s$ is completely determined by the degrees of the
components.

\begin{theorem}\label{27Dec06}
Let $K$ be a commutative $\Q$-algebra. For each $\s \in \G'$, the
signature map $\sign_\s : \Dec (\s ) \mapsto \Sign (\s)$ is a
bijection. In particular, there are only finitely many, namely
$|\Sign (\s ) |$, distinct decompositions $\s = \s_1\cdots \s_s$.
\end{theorem}

For each natural number $n\geq 1$, the set
$$\G'_n:= \{ \s \in \G'\, | \, \s (x) - x \in \sum_{i\geq 1}
Kx^{1+ni}\}$$ is a submonoid of $\G'$, $\G'= \G'_1$, and $m| n$
($m$ divides $n$) implies $\G'_n\subseteq \G'_m$.

\begin{theorem}\label{26Dec06}
Let $K$ be a commutative $\Q$-algebra. If $\s \in \G'_n$ and $\s=
\s_1\cdots \s_s$ in $\G'$ then all $\s_i\in \G'_n$.
\end{theorem}

Suppose that a monomorphism $\G'\ni \d : x\mapsto x+c_2x^2+\cdots
+ c_dx^d$, $c_i\in K$, $c_d\in K^*$, is a product $\d = \s \tau$
of two monomorphisms $\s , \tau \in \G'\backslash \{ e\}$, we say
that $\d$ is {\em decomposable}. Theorem \ref{2Dec7} gives the
{\em formulae} for $\s$ and $\tau$ via the constants $c_2, \ldots
, c_d$ of the polynomial $\d (x)$ and the degree $\deg (\s )$ of
$\s$.

A decomposability criterion for elements of $\G'$ is given
(Corollary \ref{b31Dec06}). Using it and the chain rule, necessary
conditions for irreducibility of a polynomial is found (Corollary
\ref{c2Dec7}).


\section{Proof of Theorem \ref{27Dec06}}\label{}

In this section, Theorems \ref{27Dec06} and \ref{26Dec06} are
proved based on Theorem \ref{e28Dec06}. A criterion of
decomposability (Corollary \ref{b31Dec06}) for a monomorphism of
$\G'$ is given.

A polynomial $p\in K[x]$ of the type $x+\l_2x^2+\cdots +\l_nx^n$,
$\l_i\in K$, $\l_n\in K^*$,  is called a {\em unitary} polynomial.
Note that the map $p\mapsto (\d_p : x\mapsto p)$ is a bijection
between the set of unitary polynomials and $\G'$.

If $K$ is a {\em field} of characteristic zero then the fact that
for each monomorphism $e\neq \s \in \G'$ there are only {\em
finitely many} distinct  decompositions $\s = \s_1\cdots \s_s$
follows directly from the chain rule and the fact that the
polynomial algebra $K[x]$ is a unique factorization domain.

\begin{lemma}\label{29Dec06}
Let $K$ be  a field of characteristic zero and $e\neq \s \in \G'$.
There are only  finitely many distinct  decompositions $\s =
\s_1\cdots \s_s$ in $\G'$ with all $\s_i\neq e$.
\end{lemma}

{\it Proof}. Using induction on the degree $\deg (\s )$ and
(\ref{degm}), it suffices to show that there are finitely many
distinct decompositions for $s=2$, i.e. $\s = \s_1\s_2$. Let $F:=
\s (x)$, $g:= \s_1(x)$, and $f:= \s_2(x)$. Then $F= f(g)$, the
composition of polynomials. By the chain rule, $F'= f'(g) g'$, the
product of unitary polynomials,  $g'$ is a divisor of the
polynomial $F'$ where $a':= \frac{da}{dx}$. Since there are only
finitely many  divisors of the polynomial $F'$ with the scalar
term 1 the result follows. $\Box $

{\bf The ratio form of $\s\in \G'$}. Let $K$ be a commutative
$\Q$-algebra. Each element $\s \in \G'$ is uniquely determined by
the polynomial $\s (x) = x( 1+\l_1 x+\cdots +\l_nx^n)$ where
$\l_i\in K$, $\l_n\in K^*$, and $\deg (\s ) = n$. For
computational reason it is convenient to write the polynomial $\s
(x)$ in the {\em ratio} form 
\begin{equation}\label{rfs}
\s (x) = x(1+(a_{n-1}x^{-(n-1)}+\cdots +a_1x^{-1}+1) \l x^n)= x(
1+ (A+1)\l x^n)
\end{equation}
where $\l := \l_n$; $a_i:= \frac{\l_{n-i}}{\l_n}$, $i=1, \ldots ,
n-1$; and $A:= \sum_{i=1}^{n-1} a_ix^{-i}$. Let $\R^{n-1}$ be the
standard $(n-1)$-dimensional vector space over the reals $\R$ with
the standard inner product $yz:= y_1z_1+\cdots + y_{n-1}z_{n-1}$.
Let $\alpha = (\alpha_1, \ldots , \alpha_{n-1})\in \N^{n-1}$,
$|\alpha | : =\alpha_1+\cdots +\alpha_{n-1}$, $a^\alpha :=
a_1^{\alpha_1} \cdots a_{n-1}^{\alpha_{n-1}}$, $\vn1 := (1,2,
\ldots , n-1)$, $\vn1 \alpha = \alpha_1+2\alpha_2+\cdots +
(n-1)\alpha_{n-1}$ (the inner product of vectors). For $\alpha \in
\N^{n-1}$ with $|\alpha |\leq i$, ${i\choose \alpha}:=
\frac{i!}{(i-|\alpha |)\alpha_1!\cdots \alpha_{n-1}!}$ is the
multi-binomial coefficient. We set $0^0:=1$.

For each $m\geq 1$, a direct computation gives the equality
\begin{equation}\label{sxmf}
\s (x^m) = x^m ( 1+\sum_{i=1}^m{m\choose
i}\l^ix^{in}+\sum_{i=1}^m\sum_{1\leq |\alpha | \leq i}{m\choose i}
{i\choose \alpha}a^\alpha \l^i x^{-\vn1 \alpha +in}).
\end{equation}
In more detail,
\begin{eqnarray*}
\s (x^m) &= & x^m ( 1+ (A+1)\l x^n)^m= x^m(1+\sum_{i=1}^m{m\choose i} (A+1)^i\l^ix^{in}) \\
 &=& x^m(1+\sum_{i=1}^m{m\choose i}(1+\sum_{1\leq |\alpha | \leq
 i} {i\choose \alpha} a^{\alpha } x^{-\vn1 \alpha})\l^i x^{in}),
\end{eqnarray*}
multiplying out we obtain (\ref{sxmf}).

\begin{theorem}\label{e28Dec06}
Let $K$ be a commutative $\Q$-algebra, $\G'\ni \d : x\mapsto \d
(x) =x+c_2x^2+\cdots + c_dx^d$, $c_i\in K$, $c_d\in K^*$. If $\d =
\s \tau$ for some elements $\s , \tau \in \G'$ then the elements
$\s $ and $\tau $ are uniquely determined by the degree of $\s$.
In more detail, if $\s (x) = x( 1+(a_{n-1}x^{-(n-1)}+\cdots +
a_1x^{-1}+1)\l x^n)$ where $n+1= \deg (\s)$, and $\tau (x) =
x+\mu_2 x^2+\cdots + \mu_mx^m$ where $\deg (\tau ) = m
=\frac{d}{n+1}$ then
\begin{enumerate}
\item the coefficients $a_1, \ldots , a_{n-1}$, $\l$ are the
unique solution to the following triangular system of equations:
\begin{eqnarray*}
 ma_j+\sum_{C_j}{m\choose \alpha_1, \ldots , \alpha_{j-1}}a_1^{\alpha_1}\cdots \alpha_{j-1}^{\alpha_{j-1}}
 &=& \frac{c_{d-j}}{c_d}, \;\; j = 1, \ldots , n-1,  \\
 m\l^{-1} +\sum_{\alpha \in \N^{n-1},\, |\alpha | \leq m, \, \vn1 \alpha = n} {m\choose \alpha } a^\alpha & =&
 \frac{c_{d-n}}{c_d},
\end{eqnarray*}
where $C_1:= \emptyset$ and $C_j:= \{ \alpha =(\alpha_1, \ldots ,
\alpha_{j-1})\in \N^{j-1}\, | \, \alpha_1+2\alpha_2+\cdots +
(j-1)\alpha_{j-1} = j, \, |\alpha | \leq m \}$, $j\geq 2$, and
\item the coefficients $\mu_j$, $j=2, \ldots , m$, are the unique
solution to the following triangular system of liner equations in
$\mu_j$:
$$ \mu_j +\sum_{D_j}\mu_{m'}{m'\choose
i}\l^i+\sum_{E_j}\mu_{m'}{m'\choose i} {i\choose \alpha } a^\alpha
\l^i= c_j, \;\; j=2, \ldots , m,$$ where $\mu_1:=1$,
\begin{eqnarray*}
 D_j&:= & \{ (m', i)\in \N^2\, | \, 1\leq i\leq m'<j, \; m'+in=j\}, \\
 E_j&:=& \{ (m', i, \alpha )\in \N^2\times \N^{n-1}\, | \, 1\leq
 |\alpha | \leq i\leq
m'<j,\;  m'-\vn1 \alpha +in=j\}.
\end{eqnarray*}
\end{enumerate}
\end{theorem}

{\it Proof}. Note that $\tau (x) = b+\mu_mx^m$ where $b:=
x+\sum_{i=1}^{m-1}\mu_ix^i$, $\deg \s (b) = \deg \s (x^{m-1}) =
(m-1) (n+1)$ and 
\begin{equation}\label{dgdi}
\deg \s (\mu_mx^m) -\deg \s (b) = m(n+1) - (m-1)(n+1) = n+1.
\end{equation}
 By (\ref{sxmf}) and (\ref{dgdi}), $c_d=\mu_m \l^m$, and for each
 $j=1, \ldots , n-1$,
 \begin{eqnarray*}
  c_{d-j} &=& \mu_m\l^m ( ma_j+\sum_{C_j}{m\choose \alpha_1, \ldots , \alpha_{j-1}}a_1^{\alpha_1}\cdots a_{j-1}^{\alpha_{j-1}}), \\
c_{d-n}&=& \mu_m\l^m(m\l^{-1} +\sum_{\alpha \in \N^{n-1},\,
|\alpha | \leq m, \,\vn1 \alpha =n}{m\choose \alpha }a^\alpha ).
\end{eqnarray*}
Taking the ratios $\frac{c_{d-j}}{c_d}$, $j=1, \ldots , n$, we
have the system of equations as in the first statement. The
triangular structure of the system gives the unique solution for
coefficients $a_1, \ldots , a_{n-1}, \l$. Since $\s$ is an algebra
monomorphism of $K[x]$, the element $\tau$ in the equality $\d =
\s \tau$ is unique. The second statement gives the exact values
for the coefficients by applying the Cramer's formula.

2. A direct computation gives
\begin{eqnarray*}
 \d (x) &=&\sum_{m'=1}^m\mu_{m'}x^{m'}+\sum_{m'=1}^m\sum_{i=1}^{m'}\mu_{m'}{m'\choose i} \l^ix^{m'+in}\\
 &+&
 \sum_{m'=1}^m\sum_{i=1}^{m'}\sum_{1\leq |\alpha | \leq i}\mu_{m'}{m'\choose i}{i\choose \alpha } a^\alpha \l^ix^{m'-\vn1 \alpha+in}.
\end{eqnarray*}
In more detail,
\begin{eqnarray*}
 \d(x)&=&\sum_{m'=1}^m\mu_{m'}\s (x^{m'})= \sum_{m'=1}^m\mu_{m'}x^{m'} (1+\sum_{i=1}^{m'}{m'\choose i}
 ( 1+ \sum_{1\leq |\alpha |\leq i} {i\choose \alpha } a^\alpha x^{-\vn1 \alpha } ) \l^ix^{in}). \\
 &
\end{eqnarray*}
By opening up the brackets we obtain the result. The formula for
$\d (x)$ above can be written shortly as the sum $S_1+S_2+S_3$
where $S_1$, $S_2$ and $S_3$ are the single, double, and triple
sum respectively. For each $j=2, \ldots , m$, the coefficient of
$x^j$ in the sum $S_1$ is equal to $\mu_j$; the coefficient of
$x_j$ in the sum $S_2$ is equal to $\sum_{D_j} \mu_{m'}{m'\choose
i}\l^i$ since the conditions $1\leq i\leq m'$ and $m'+in=j$ imply
$m'<j$; and the coefficient of $x^j$ in the sum $S_3$ is equal to
$\sum_{E_j} \mu_{m'} {m'\choose i} {i\choose \alpha} a^\alpha
\l^i$ since the two conditions $in-\vn1 \alpha \geq 1$ (as follows
easily from the equality $\alpha_1+\cdots + \alpha_{n-1} \leq i$)
and $ m'-\vn1 \alpha +in=j$ imply the strict inequality $m'<j$. In
more detail,
\begin{eqnarray*}
 in-\vn1 \alpha &=&in-\alpha_1-2\alpha_2-\cdots - (n-1)\alpha_{n-1} \\
 &\geq & n(\alpha_1+\cdots +\alpha_{n-1})-\alpha_1-2\alpha_2-\cdots - (n-1)\alpha_{n-1} \\
 &=&\sum_{i=1}^{n-1}(n-i)\alpha_i\geq 1
\end{eqnarray*}
since $|\alpha | \geq 1$ and all $n-i\geq 1$. For each $j=2,
\ldots , m$, equating the coefficients of $x^j$ of both sides of
the equality $\s \tau (x) = \d (x)$ we obtain the triangular
system of linear equations with the unknowns $\mu_j$ as in
statement 2. It has the unique solution that can be easily written
explicitly using the Cramer's formula via the elements $a_1,
\ldots , a_{n-1}, \l$. This proves the theorem. $\Box $

Equating the coefficients of $x^j$, $m<j<d-n$, of both sides of
the equality $\d (x) = \s\tau (x)$, we get the system of equations
\begin{equation}\label{cjrem}
c_j=\sum_{F_j}\mu_{m'}{m'\choose i}
\l^i+\sum_{G_j}\mu_{m'}{m'\choose i}{i\choose \alpha} a^\alpha
\l^i,\;\;\; m<j<d-n,
\end{equation}
where
\begin{eqnarray*}
 F_j&:=&\{ (m',i)\in \N^2\, | \, 1\leq i\leq m'\leq m, \;
m'+in=j\}, \\
 G_j&:=&\{ (m',i, \alpha )\in \N^2\times \N^{n-1}\, |
\, 1\leq |\alpha |\leq i\leq  m'\leq m, \;  m'-\vn1 \alpha
+in=j\}.
\end{eqnarray*}
Note that in (\ref{cjrem}),
$$ a^\alpha \l^i =
\l_{n-1}^{\alpha_1}\l_{n-2}^{\alpha_2}\cdots
\l_1^{\alpha_{n-1}}\l^{i-\alpha_1-\cdots - \alpha_{n-1}}, \;\;
i-\alpha_1-\cdots - \alpha_{n-1}\geq 0.$$ So, the RHS of
(\ref{cjrem}) is a polynomial    in $\mu_i$ and $\l_k$ with
integer coefficients.

{\bf  Proof of Theorem \ref{27Dec06}}. By the very definition, the
signature map is surjective. If $ \s = \s_1\cdots \s_s=
\tau_1\cdots \tau_s$ in $\G'$ and $\deg (\s_1) = \deg (\tau_1),
\ldots , \deg (\s_s) =\deg (\tau_s)$ then, by Theorem
\ref{e28Dec06}, $\s_1=\tau_1, \ldots , \s_s=\tau_s$,  i.e. the
signature map is injective. $\Box$

\begin{corollary}\label{a31Dec06}
We keep the notation of Theorem \ref{e28Dec06}. Then
\begin{enumerate}
\item for each $j=1, \ldots , n-1$, $a_j\in \Z
[\frac{1}{m!}][\frac{c_{d-j}}{c_d}, \frac{c_{d-j+1}}{c_d}, \ldots
, \frac{c_{d-1}}{c_d}]$, i.e. $a_j$ is a polynomial in
$\frac{c_{d-j}}{c_d}, \frac{c_{d-j+1}}{c_d}, \ldots ,
\frac{c_{d-1}}{c_d}$ with coefficients from $\Z [\frac{1}{m!}]$.
\item $\l \in \Q (\frac{c_{d-n}}{c_d}, \frac{c_{d-n+1}}{c_d},
\ldots , \frac{c_{d-1}}{c_d})$, i.e. $\l $ is a rational function
in $\frac{c_{d-n}}{c_d}, \frac{c_{d-n+1}}{c_d}, \ldots ,
\frac{c_{d-1}}{c_d}$ with rational coefficients, \item  $\mu_j\in
\sum_{k=2}^jc_k\Q (\frac{c_{d-n}}{c_d}, \frac{c_{d-n+1}}{c_d},
\ldots , \frac{c_{d-1}}{c_d})$ and $\mu_j\in \sum_{k=2}^jc_k\Q
(\l_1,\ldots , \l_n )$, for each $j=2, \ldots , m$.
\end{enumerate}
\end{corollary}

{\it Proof}. It is obvious.  $\Box $

Recall  that a monomorphism $\d\in \G'\backslash \{ e\}$ is {\em
decomposable} if $\d = \s \tau$ for some monomorphisms $\s , \tau
\in \G'\backslash \{ e\}$.  The next corollary is a
decomposability criterion.

\begin{corollary}\label{b31Dec06}
We keep the notation of Theorem \ref{e28Dec06}. The monomorphism
$\d$ of $\G'$ of degree $d=(n+1)m$, $\d (x) = x+c_2x^2+\cdots
+c_dx^d$ is equal to the product $\d = \s \tau$ of some
monomorphisms $\s$ and $\tau $ (as in Theorem \ref{e28Dec06}) with
$\deg (\s ) = n+1$ and $\deg (\tau ) = m$ iff $c_d=\mu_m\l^m$ and
the coefficients $c_j$, $ m<j<d-n$, satisfy the equation
(\ref{cjrem}) (i.e. the coefficients $c_d$ and $c_j$ are uniquely
determined by the elements $c_2, \ldots , c_m,
\frac{c_{d-n}}{c_d}, \frac{c_{d-n+1}}{c_d}, \ldots ,
\frac{c_{d-1}}{c_d}$).
\end{corollary}

{\it Proof}. In the proof of Theorem \ref{e28Dec06} we used only
the equalities of the coefficients of $x^j$ in $\d (x) = \s \tau
(x)$ for $j=2, \ldots , m, d-n, d-n+1, \ldots , d-1$. The
remaining equalities are $c_d= \mu_m\l^m$ for $j=d$ and
(\ref{cjrem}) for $j$ such that $m<j<d-n$. Now, the corollary is
obvious.  $\Box $

{\bf Proof of Theorem \ref{26Dec06}}. We prove the theorem in two
steps: first,  when the algebra $K$ contains a primitive $n$'th
root of unity, and then the general case can be reduced to the
first one using Corollary \ref{a31Dec06}.

 Suppose that the algebra $K$ contains a primitive $n$'th
root of unity, say $\l$.  $\Aut_K(K[x])\ni \g : x\mapsto \l x$,
and $\o= \o_\g : \tau \mapsto \g \tau \g^{-1}$ is  the inner
automorphism of the group $\Aut_{K, c}(K[[x]])$ of {\em
continuous} automorphisms of the series algebra $K[[x]]$. It is
obvious that $\G'_n= \G'^\o:= \{ \d \in \G' \, | \, \o (\d ) = \d
\}$. Since
$$ \s_1\cdots \s_s=\s = \o ( \s ) = \o (\s_1) \cdots \o (\s_s)$$
and $\deg (\s_1) = \deg (\o (\s_1)), \ldots , \deg (\s_s) = \deg
(\o (\s_s))$ we must have $\s_1=\o(\s_1), \ldots , \s_s= \o
(\s_s)$, by Theorem \ref{27Dec06}, i.e. $\s_1, \ldots , \s_s\in
\G'^\o = \G'_n$, as required.

In the general case, fix a commutative $\Q$-algebra, say $L$, that
contains both $K$ and a primitive $n$'th root of unity. Then
$\G'\subseteq \G'(L)$ and $\G'_n\subseteq \G'_n(L)$. By the
previous case, $\s_1, \ldots , \s_s\in \G'_n(L)$, and, by
Corollary \ref{a31Dec06}, $\s_1, \ldots , \s_s\in \G'_n$.
Therefore, $\s_1, \ldots , \s_s\in \G'\cap \G'_n(L)=\G'_n$, as
required. $\Box$

Let $\Aut_{\Q}(K)$ be the group of $\Q$-algebra automorphisms of
$K$. Each element $\th\in\Aut_{\Q}(K)$ acts naturally on the
polynomial algebra $K[x]$, $\th (\sum_{i\geq
0}\l_ix_i)=\sum_{i\geq 0} \th (\l_i) x^i$. Let $\Aut (\G')$ be the
 group of automorphisms of the monoid $\G'$. The map
 $$ \o_\cdot :\Aut_{\Q}(K)\ra \Aut (\G'), \;\; \th\mapsto \o_\th : \s \mapsto \th \s \th^{-1},$$
is a group monomorphism as follows from the equality $(\th
\s\th^{-1}) (x) = \th ( \s (x))$.

\begin{corollary}\label{t6Jan7}
Let $K$ be a commutative $\Q$-algebra, $\s \in \G'$ and $\th \in
\Aut_{\Q}(K)$. If $ \s = \s_1\cdots \s_s$ in $\G'$ and $\o_\th (\s
)=\s$ then $\o_\th (\s_1 )=\s_1, \ldots , \o_\th (\s_ s)=\s_s$.
\end{corollary}

{\it Proof}. Since $\deg \, \o_\th (\tau ) = \deg \, \tau$ for all
$\tau \in \G'$ and
$$ \s_1\cdots \s_s= \s = \o_\th (\s ) = \o_\th (\s_1) \cdots
\o_\th (\s_s)$$ we must have $\o_\th (\s_1 )=\s_1, \ldots , \o_\th
(\s_ s)=\s_s$, by Theorem \ref{e28Dec06}. $\Box $


\section{Formulae for the components $\s$ and $\tau$ in $\d = \s\tau$}\label{}

In this section, $K$ is a commutative $\Q$-algebra. If a
monomorphism $\d \in \G'$ is decomposable, $\d = \s\tau$, then one
can write formulae for the monomorphisms $\s$ and $\tau$ via the
coefficients of the polynomial $\d (x)$ and the degree $\deg (\s
)$ of $\s$ (Theorem \ref{2Dec7}). In order to do so, we use the
inversion formula \cite{Bav-inform} for an automorphism of a
polynomial algebra $P_n:= K[x_1, \ldots , x_n]$.

{\bf The inversion formula}. Let us recall the inversion formula
(for details the reader is referred to \cite{Bav-inform}). For
purpose of application (because we have $n-1$ elements $a_1,
\ldots , a_{n-1}$ in Theorem \ref{e28Dec06}), it is convenient to
state it for a polynomial algebra $P_{n-1}:=K[x_1, \ldots ,
x_{n-1}]$ in $n-1$ variables rather than $n$ as in
\cite{Bav-inform}.

An automorphism $s\in \Aut_K(P_{n-1})$ is uniquely determined by
the elements $x_1':= s(x_1), \ldots , x_{n-1}':=s(x_{n-1})$. Let
$\der_1:=\frac{\der}{\der x_1}, \ldots ,
\der_{n-1}:=\frac{\der}{\der x_{n-1}}$ be the partial derivatives
that corresponds to the canonical generators $x_1, \ldots ,
x_{n-1}$ of the polynomial algebra $P_{n-1}$. The corresponding to
the elements $x_1',\ldots , x_{n-1}'$  partial derivatives
$\der_1':=\frac{\der}{\der x_1'}, \ldots ,
\der_{n-1}':=\frac{\der}{\der x_{n-1}'}$ are given by the rule
\begin{equation}\label{MPdad2}
\der_i' (\cdot ) := \D^{-1} \det
 \begin{pmatrix}
  \frac{\der s (x_1)}{\der x_1} & \cdots & \frac{\der s (x_1)}{\der x_{n-1}} \\
  \vdots & \vdots & \vdots \\
\frac{\der }{\der x_1} (\cdot )& \cdots & \frac{\der }{\der x_m}(\cdot )\\
 \vdots & \vdots & \vdots \\
\frac{\der s (x_{n-1})}{\der x_1} & \cdots & \frac{\der s (x_{n-1})}{\der x_{n-1}} \\
\end{pmatrix}, \;\;\; i=1, \ldots , n-1,
\end{equation}
where we `drop' $s (x_i)$ in the {\em Jacobian}  $\D :=\det
(\frac{\der s (x_i)}{\der x_j})\in K^*$.

For each $i=1, \ldots , n-1$, and $j\geq 0$,  let
\begin{equation}\label{MhPdad3}
\phi_i':= \sum_{k\geq 0}(-1)^k\frac{x_i'^k}{k!}\der_i'^k, \;\;\;
\phi_{i,j}':= \sum_{k=0}^j(-1)^k\frac{x_i'^k}{k!}\der_i'^k:
P_{n-1} \ra P_{n-1},
\end{equation}
and 
\begin{equation}\label{MhPdad4}
\phi_s := \phi_1'\cdots \phi_{n-1}':P_{n-1} \ra P_{n-1}, \;\;
\phi_s(P_{n-1})=K.
\end{equation}

\begin{theorem}\label{MhPi8Nov05}
{\rm (The Inversion Formula, \cite{Bav-inform})} For each $s \in
\Aut_K(P_{n-1})$
 and $a\in P_{n-1}$,
 $$ s^{-1}(a)=\sum_{\alpha \in \mathbb{N}^{n-1}}\phi_s
 (\frac{\der'^\alpha}{\alpha!}(a))x^\alpha  $$
 where $\phi_s
 (\frac{\der'^\alpha}{\alpha!}(a))\in K$, $x^\alpha := x_1^{\alpha_1}\cdots x_{n-1}^{\alpha_{n-1}}$
 and $\der'^\alpha := {\der'}_1^{\alpha_1}\cdots
 {\der'}_{n-1}^{\alpha_{n-1}}$.
\end{theorem}

{\it Remark}. In \cite{Bav-inform}, the theorem is proved for a
field $K$ of characteristic zero but the proof goes without a
change for an arbitrary commutative $\Q$-algebra $K$.

{\bf The formulae for $\s$ and $\tau$ in $\d =\s\tau$}. Let $\d =
\s\tau$ be as in Theorem \ref{e28Dec06}. Using the inversion
formula (Theorem \ref{MhPi8Nov05}) we obtain the formulae for the
monomorphisms $\s$ and $\tau $ (Theorem \ref{2Dec7}).

Consider the automorphism $s: P_{n-1}\ra P_{n-1}$ given by the
expressions for the elements  $a_j$ in Theorem \ref{e28Dec06}:
$$x_j':= s(x_j) := mx_j+\sum_{C_j}{m\choose \alpha_1, \ldots ,
\alpha_{n-1}}x_1^{\alpha_1} \cdots x_{j-1}^{\alpha_{j-1}}, \;\;
j=1, \ldots , n-1.$$ Its triangular structure guarantees that $s$
is an automorphism of $P_{n-1}$ with the Jacobian $\D := \det
(\frac{\der s(x_i)}{\der x_j})=m^{n-1}$. Putting $a=x_j$ in
Theorem \ref{MhPi8Nov05} and then applying $s$ yields
\begin{equation}\label{xjf1}
x_j=\sum_{\alpha \in \N^j}\phi_s(\frac{\der'^\alpha}{\alpha
!}(x_j))x'^\alpha
\end{equation}
where $x'^\alpha := s (x^\alpha  )= x_1'^{\alpha_1}\cdots
x_j'^{\alpha_j}$ (we used also the triangular structure of $s$).
Each $x_j$ is a polynomial in variables $x_1', \ldots , x_j'$ with
coefficients $ \phi_s(\frac{\der'^\alpha}{\alpha !}(x_j))\in K$.
Therefore, instead of the map $\phi_s$ in (\ref{xjf1}) one can
write the map $\phi'_1\cdots \phi_j'$, i.e. 
\begin{equation}\label{xjf}
x_j=\sum_{\alpha \in \N^j}\phi_1'\cdots
\phi_j'(\frac{\der'^\alpha}{\alpha !}(x_j))x'^\alpha .
\end{equation}

The {\em total} degree $\deg_{x'}(x_j)$ of the polynomial $x_j$
with respect to the variables $x_1', \ldots , x_{n-1}'$ (or  to
$x_1', \ldots , x_j'$) satisfies the inequality 
\begin{equation}\label{dxpi}
\deg_{x'}(x_j)\leq j.
\end{equation}
To prove this inequality we use induction on $j$. The case $j=1$
is obvious since $ x_1'= m x_1$. Suppose that $j\geq 2$ and the
inequality is true for all $j'<j$. Since $mx_j =
x_j'-\sum_{C_j}{m\choose \alpha_1, \ldots ,
\alpha_{j-1}}x_1^{\alpha_1}\cdots x_{j-1}^{\alpha_{j-1}}$ we have
$$ \deg_{x'}(x_j) \leq \sum_{k=1}^{j-1} \alpha_k
\deg_{x'}(x_k)=\sum_{k=1}^{j-1}\alpha_k k \leq j,$$ see the
definition of the set $C_j$ in Theorem \ref{e28Dec06}. By
induction on $j$, (\ref{dxpi}) holds.

Note that by (\ref{xjf}),
$$ \deg_{x'}(\frac{\der'^\alpha}{\alpha !}(x_j))\leq
\deg_{x'}(x_j) - |\alpha | \leq j-|\alpha |.
$$

In a view of (\ref{dxpi}) and the triangular structure of the
automorphism $s$, the formula for $x_j$, (\ref{xjf}), can be
written as follows 
\begin{equation}\label{xjf2}
x_j=x_j(x_1', \ldots , x_j')=\sum_{\alpha \in \N^j, |\alpha |\leq
j}\phi_{1,j-|\alpha |}'\cdots \phi_{j,j-|\alpha |}'
(\frac{\der'^\alpha}{\alpha !}(x_j))x'^\alpha ,
\end{equation}
it contains only finitely many terms where $x'^\alpha :=
x_1'^{\alpha_1}\cdots x_j'^{\alpha_j}$.

\begin{theorem}\label{2Dec7}
We keep the notation of Theorem \ref{e28Dec06}. Then
\begin{enumerate}
\item the coefficients $\l$ and $a_j$, $j=1, \ldots , n-1$, are
given by the rule
\begin{eqnarray*}
a_j &=& x_j(\frac{c_{d-1}}{c_d}, \ldots , \frac{c_{d-j)}}{c_d})
=\sum_{\alpha \in \N^j, |\alpha |\leq j}\phi_{1,j-|\alpha
|}'\cdots \phi_{j,j-|\alpha |}' (\frac{\der'^\alpha}{\alpha
!}(x_j))(\frac{c_{d-1}}{c_d})^{\alpha_1} \cdots  (\frac{c_{d-j}}{c_d})^{\alpha_j} ,\\
 \l &=& m^{-1}( \frac{c_{d-n}}{c_d}-\sum_{\alpha \in \N^{n-1},\, |\alpha | \leq m, \,  \vn1 \alpha =n} {m\choose \alpha
 }(\frac{c_{d-1}}{c_d})^{\alpha_1} \cdots
 (\frac{c_{d-(n-1)}}{c_d})^{\alpha_{n-1}})^{-1}.
\end{eqnarray*}
\item Then the coefficients $\mu_j$, $j=2, \ldots , m$, can be
written explicitly via the elements $a_1, \ldots , a_{n-1}, \l$
using the Cramer's formula for the unique solution of a system of
linear equations.
\end{enumerate}
\end{theorem}

{\it Proof}. 1. This follows at once from (\ref{xjf2}) and the
formulae for the elements $a_j$ and $\l $ in Theorem
\ref{e28Dec06}.

2. It is obvious.  $\Box $

\begin{theorem}\label{23Dec06}
Let $K$ be a commutative $\Q$-algebra. Then the submonoid of
$\G'$, say $M$,  generated by the set $\{ \s_{\l x^n}:x\mapsto
x+\l x^n\, | \, \l\in K^*, n\geq 2\}$ is a free monoid, i.e.
$\s_{a_1}\cdots \s_{a_s} = \s_{b_1}\cdots \s_{b_t}$ iff $s=t$,
$a_1=b_1, \ldots , a_s=b_s$.
\end{theorem}

{\it Proof}. Let $a=\l x^{n+1}$, $\l \in K^*$, $n\geq 1$. Then
$\s_a(x)= x(1+\l x^n)$ and, for each $m\geq 1$,
$$ \s_a(x^m) = x^m(1+\l x^n)^m= x^m\sum_{i=0}^m{m\choose i}
\l^ix^{in}.$$ The polynomial $\s_a(x^m)$ has degree $m+mn$. It is
 the  sum of monomials ${m\choose i} \l^ix^{m+in}$ with coefficients
from $K^*$, the {\em leading} and the {\em pre-leading} terms are
$l:= \l^mx^{m+mn}$ and $p:= m\l^{m-1} x^{m+(m-1)n}$ respectively.
Note that $\deg \, \s_a(x) <\deg\, \s_a(x^2) <\cdots < \deg \,
\s_a(x^i)<\cdots $ and
$$ \deg \, \s_a(x^{m-1})= m-1+(m-1)n<m+(m-1)n= \deg \, p.$$
Therefore, for any polynomial, say $f=\mu x^m+\cdots $,  of degree
$m$ with the leading coefficient $\mu \in K^*$, $\s_a(f) =
L+P+\cdots $ where $L:= \mu l$ and $P:=\mu p$ are the leading and
the pre-leading terms of the polynomial $\s_a(f)$ and the three
dots mean smaller terms. Note that $\frac{L}{P}=\frac{\l}{m}x^n$,
hence 
\begin{equation}\label{nml}
n=\deg (\frac{L}{P}), \;\; m = \frac{\deg \, \s_a(f)}{n+1}, \;\;
\l = mx^{-n} \frac{L}{P}.
\end{equation}
Let $\s= \s_{a_1} \cdots \s_{a_s}$. By induction on $s$, it is
easy to prove that the coefficients of the leading and the
pre-leading term of the element $\s (x)$ are units of $K$. By
(\ref{nml}), the element $a_1$ is uniquely determined by the
leading and the pre-leading terms of the polynomial $\s (x)$.
Since $\s = \s_{b_1}\cdots \s_{b_t}$, we must have $a_1=b_1$. Then
$\s_{a_1} = \s_{b_1}$ and then $\s_{a_2} \cdots
\s_{a_s}=\s_{b_2}\cdots \s_{b_t}$. Repeating the same argument we
see that the equality $\s_{a_1} \cdots \s_{a_s}=\s_{b_1}\cdots
\s_{b_t}$ holds iff $s=t$, $a_1=b_1, \ldots , a_s=b_s$. $\Box $


\section{Necessary conditions for irreducibility of a polynomial}\label{}

In this section, necessary conditions for irreducibility of a
polynomial is given (Corollary \ref{c2Dec7}). These conditions are
far from being sufficient. Their advantage is that they are
explicit and if they do not hold then one can find explicitly
(using Theorem \ref{2Dec7}) a nontrivial divisor of the polynomial
(i.e. a formula for certain divisors is given explicitly).

Let $K$ be a commutative $\Q$-algebra. Let $\CP$ be a set of all
the polynomials of the type $p=1+\l_1x+\cdots +\l_tx^t$ where
$\l_i\in K$, $\l_t\in K^*$, and $\CQ$ be the set of all the
unitary polynomials, i.e. of the type $q= x+\mu_2x^2+\cdots +
\mu_sx^s$ where $\mu_i\in K$ and $\mu_s\in K^*$. The derivation
$(\cdot )':= \frac{d}{dx}:\CQ\ra \CP$, $q\mapsto q'$, is a {\em
bijection} with the inverse
$$\int :\CP\ra \CQ , \;\; 1+\sum_{i\geq 1}\l_i x_i\mapsto
x+\sum_{i\geq 1} \frac{\l_i}{i+1}x^{i+1}.$$

A polynomial $p\in \CP$ is called {\em reducible} if $p=qr$ for
some polynomial $q,r\in \CP$ each of degree $\geq 1$. A polynomial
$p\in \CP$ of degree $d-1\geq 1$ can be written uniquely in the
form $p= 1+\sum_{i=2}^d \frac{c_i}{i-1}x^{i-1}$. Then its integral
$\int p := x+c_2x^2+\cdots + c_dx^d$ determines the monomorphism
$\G'\ni \d_p : x\mapsto \d_p(x) = \int p$. We say that a
polynomial $p\in \CP$ is {\em decomposable} if the monomorphism
$\d_p$ is {\em decomposable}, i.e. $\d_p= \s\tau$ for some
monomorphisms $\s , \tau \in \G'\backslash \{ e\}$. It is obvious
that each polynomial of prime degree is indecomposable, i.e. not
decomposable.

Let $f(x):= \s (x)$ and $g(x) := \tau (x)$. The polynomials $f$
and $g$ are {\em unitary} and have degree $\geq 2$. Therefore,
their derivatives $f'$ and $g'$ belong to $ \CP$ and have degree
$\geq 1$. By the chain rule, 
\begin{equation}\label{pchr}
p:=\d(x)'=(f(g(x)))'= f'(g(x))\cdot g'.
\end{equation}
\begin{lemma}\label{i2Dec7}
\begin{enumerate}
\item A polynomial $p\in \CP$ is decomposable iff $p= f'(g)g'$ for
some  polynomials $f,g\in \CQ$ of degree $\geq 2$. \item Each
decomposable polynomial is reducible.
\end{enumerate}
\end{lemma}

{\it Proof}. Both statements follow at once from (\ref{pchr}).
$\Box $

{\it Remarks}.  1. If the polynomial $p\in \CP$ is decomposable
then using Theorem \ref{2Dec7} one can find explicitly all the
pairs $(f'(g), g')$ of the divisors of $p$ as in Lemma
\ref{i2Dec7}, i.e. $p=f'(g)g'$.

2. Moreover, using Theorem \ref{2Dec7} and Corollary
\ref{b31Dec06}, one can find the set of all decomposable
polynomials explicitly.

 The next corollary gives necessary
conditions for a polynomial being irreducible.

\begin{corollary}\label{c2Dec7}
Let $K$ be a commutative $\Q$-algebra. If a polynomial $p=
1+\sum_{i=2}^d\frac{c_i}{i-1}x^{i-1}$, $c_i\in K$, $c_d\in K^*$,
is irreducible  then for the polynomial $\d (x) := \int p =
x+c_2x^2+\cdots + c_dx^d$ the conditions of Corollary
\ref{b31Dec06} do not hold for each pair $(m,n)$ such that $m\geq
2$, $n\geq 1$, $ d=m(n+1)$.
\end{corollary}

{\it Proof}. This follows directly from Lemma \ref{i2Dec7} and
Corollary \ref{b31Dec06}.  $\Box $

Department of Pure Mathematics

University of Sheffield

Hicks Building

Sheffield S3 7RH

UK

email: v.bavula@sheffield.ac.uk

\end{document}